\documentclass[conference]{IEEEtran}
\IEEEoverridecommandlockouts
\usepackage{graphicx}
\usepackage{epstopdf}
\graphicspath{{../}}
\usepackage{multirow}
\usepackage{float}
\usepackage{amsmath,amsthm}
\usepackage{amssymb}
\usepackage{array}
\usepackage{mdwmath}
\usepackage{mdwtab}
\usepackage{eqparbox}
\usepackage{url}
\usepackage{nomencl}
\usepackage{cite}
\usepackage{algorithm}
\usepackage{algorithmic}
\usepackage{makeidx}
\usepackage{ifthen}
\usepackage{subfigure}
\usepackage{caption}
\usepackage{pdflscape}
\usepackage{hyperref}
\usepackage{xcolor}

\def\mbi{\mathbb{I}}

\def\Nikin{\mathcal{N}_{ik}^{\mathrm{in}}}
\def\Nikout{\mathcal{N}_{ik}^{\mathrm{out}}}

\def\a{{\alpha}}

\def\1{\mathbf{1}}
\def\0{\mathbf{0}}
\def\la{\langle}
\def\ra{\rangle}

\newcommand{\SABTV}{$\mathcal{S\!AB}$--\textnormal{TV}}

\usepackage{amsmath,mathtools, amssymb, bbm, xspace}
\usepackage{booktabs,ragged2e}
\usepackage{multirow}
\usepackage[flushleft]{threeparttable}

\usepackage{graphicx}

\def\exp{\mathop{\hbox{\rm exp}}}

\def\spose#1{\hbox to 0pt{#1\hss}}

\def\text #1{\hbox{\quad#1\quad}}

\def\bbG{\mathbb{G}}

\def\nthinsp{\mskip -2   mu}

\def\F{\mathcal{F}}
\def\G{_{\scriptscriptstyle G}}

\def\superstar{^{\raise 0.5pt\hbox{$\nthinsp *$}}}
\def\SUPERSTAR{^{\raise 0.5pt\hbox{$*$}}}

\def\lamstarT {\lambda^{\raise 0.5pt\hbox{$\nthinsp *$}T}}

\def\hbar{\skew{4.2}\bar h}

		\def\bkE{{\rm I\kern-.17em E}}
		\def\bk1{{\rm 1\kern-.17em l}}
		\def\bkD{{\rm I\kern-.17em D}}
		\def\bkR{{\rm I\kern-.17em R}}
		\def\bkP{{\rm I\kern-.17em P}}
		\def\bkY{{\bf \kern-.17em Y}}
		\def\bkZ{{\bf \kern-.17em Z}}

		\def\bc{\begin{center}}
		\def\be{\begin{enumerate}}
		\def\bi{\begin{itemize}}
		\def\ec{\end{center}}
		\def\ee{\end{enumerate}}
		\def\ei{\end{itemize}}
		\def\es{\end{small}}
		\def\eS{\end{slide}}

	\def\cp2problem#1#2#3#4{\fbox
		 {\begin{tabular*}{0.9\textwidth}
			{@{}l@{\extracolsep{\fill}}l@{\extracolsep{6pt}}l@{\extracolsep{\fill}}c@{}}
				#1 & & $#4 $ 
			\end{tabular*}}}

		\def\bkE{{\rm I\kern-.17em E}}
		\def\bk1{{\rm 1\kern-.17em l}}
		\def\bkD{{\rm I\kern-.17em D}}
		\def\bkR{{\rm I\kern-.17em R}}
		\def\bkP{{\rm I\kern-.17em P}}
		
		\def\bkZ{{\bf{Z}}}

\newcommand {\beeq}[1]{\begin{equation}\label{#1}}
\newcommand {\eeeq}{\end{equation}}
\newcommand {\bea}{\begin{eqnarray}}
\newcommand {\eea}{\end{eqnarray}}

\def\texitem#1{\par\smallskip\noindent\hangindent 25pt
               \hbox to 25pt {\hss #1 ~}\ignorespaces}

\newcommand{\argmin}{\operatornamewithlimits{argmin}}

\newcommand{\beq}{\begin{equation}}
\newcommand{\eeq}{\end{equation}}
\newcommand{\beqn}{\begin{eqnarray}}
\newcommand{\eeqn}{\end{eqnarray}}
\newcommand{\beqno}{\begin{eqnarray*}}
\newcommand{\eeqno}{\end{eqnarray*}}
\newcommand{\bma}{\begin{displaymath}}
\newcommand{\ema}{\end{displaymath}}
\newcommand{\bnu}{\begin{enumerate}}
\newcommand{\enu}{\end{enumerate}}
\newcommand{\bce}{\begin{center}}
\newcommand{\ece}{\end{center}}
\newcommand{\btb}{\begin{tabular}}
\newcommand{\etb}{\end{tabular}}

\def\G{{\mathbb{G}}}

\def\ba{{\mathbf{a}}}
\def\bb{{\mathbf{b}}}

\def\bx{{\mathbf{x}}}
\def\by{{\mathbf{y}}}

\def\bu{{\mathbf{u}}}
\def\bv{{\mathbf{v}}}

\def\b1{{\mathbf{1}}}
\def\cNinik{{\mathcal{N}^{\rm in}_{ik}}}
\def\cNoutik{{\mathcal{N}^{\rm out}_{ik}}}
\def\cNini{{\mathcal{N}^{\rm in}_{i}}}
\def\cNouti{{\mathcal{N}^{\rm out}_{i}}}

\usepackage[normalem]{ulem} 

\newtheorem{theorem}{Theorem}

\newtheorem{proposition}{Proposition}

\newtheorem{remark}{Remark}

\newtheorem{assumption}{Assumption}

\newcommand{\bg}{{\mathbf{g}}}

\newcommand{\g}{{\mathbf{g}}}

\newcommand{\E}{{\mathbb{E}}}

\newcommand{\one}{{\mathbf{1}}}

\newcommand{\T}{{\mathsf{T}}}

\definecolor{myBlue}{rgb}{0.80,0.85,1.00}
\definecolor{myYellow}{rgb}{0.951,1.000,0.547}

\def\la{{\langle}}
\def\ra{{\rangle}}

\def\bit{\begin{itemize}}
\def\eit{\end{itemize}}
\def\BEAS{\begin{eqnarray*}}
\def\EEAS{\end{eqnarray*}}

\def\re{{\mathbb R}}

\def\a{\alpha}
\def\b{\beta}

\def\g{\gamma}

\def\argmin{\mathop{\rm argmin}}

\newcommand{\bxi}{{\boldsymbol{\xi}}}

\newcommand{\bdelta}{{\boldsymbol{\delta}}}

\makeatletter
\def\munderbar#1{\underline{\sbox\tw@{$#1$}\dp\tw@\z@\box\tw@}}

\makeatother

\allowdisplaybreaks
\begin{document}

\title{Distributed Stochastic Optimization with Gradient Tracking over Time-Varying Directed Networks\\
\thanks{This material is based in part upon work supported by the NSF award CCF-2106336 and the ARPA-H award SP4701-23-C-0074.}
}

\author{\IEEEauthorblockN{Duong Thuy Anh Nguyen, Duong Tung Nguyen, and Angelia Nedi\'c}
\IEEEauthorblockA{School of Electrical, Computer and Energy Engineering, Arizona State University, AZ, USA \\
Email: dtnguy52@asu.edu, duongnt@asu.edu, Angelia.Nedich@asu.edu}
}

\maketitle

\begin{abstract}
We study a distributed method called \SABTV, which employs gradient tracking to collaboratively minimize the sum of smooth and strongly-convex local cost functions for networked agents communicating over a time-varying directed graph. Each agent, assumed to have access to a stochastic first-order oracle for obtaining an unbiased estimate of the gradient of its local cost function, maintains an auxiliary variable to asymptotically track the stochastic gradient of the global cost. The optimal decision and gradient tracking are updated over time through limited information exchange with local neighbors using row- and column-stochastic weights, guaranteeing both consensus and optimality. With a sufficiently small constant step-size, we demonstrate that, in expectation, \SABTV converges linearly to a neighborhood of the optimal solution. Numerical simulations illustrate the effectiveness of the proposed algorithm.
\end{abstract}

\begin{IEEEkeywords}
Stochastic optimization, distributed optimization, gradient tracking, time-varying directed graphs
\end{IEEEkeywords}

\section{Introduction}

The emergence of the Internet of Things has connected devices and systems in unprecedented ways, and edge computing has become a game-changer, providing low-latency and ubiquitous computations by pushing computation and data storage to edge devices such as sensors, smartphones, and other smart devices \cite{duong_wfiot,ella23,cheng23}.
Edge computing enables more efficient and responsive systems by processing data closer to the source, offering an ideal infrastructure for distributed learning.
Distributed optimization in large-scale, multi-agent systems have thus gained substantial traction in various fields,
including sensor networks \cite{Rabbat2004}, crowdsourcing \cite{wiopt23},  machine learning \cite{Tsianos2012ML}, big-data analytics \cite{Daneshmand2015}, and distributed control \cite{Soomin2018}.

The presence of noisy and inexact data exchange at the edge may pose significant challenges for distributed optimization. In addition, the communication networks among agents can be dynamic, as they may experience communication delays, user mobility, and the impact of straggler effects, which can arise due to the requirement of synchronized communication. To address these challenges, this paper aim to investigate distributed optimization schemes that can handle inexact data and directed time-varying communication networks.
Specifically, we propose  \SABTV, a distributed stochastic optimization algorithm that minimizes the sum of local cost functions of agents communicating over a directed time-varying graph. This algorithm assumes access to a stochastic first-order oracle to obtain unbiased gradient estimates of agents' local cost functions. 
The optimal solution estimates are exchanged over time-varying row-stochastic weight matrices, while the gradient tracking are exchanged over time-varying column-stochastic weight matrices. 

The literature on deterministic finite-sum problems is extensive and includes several works such as \cite{Nedic2009, Nedic2011, Ram2009, Srivastava2011, Duchi2012, Wei2014}, while stochastic problems have also been studied \cite{Ram2009_2,Nedic2014,Ram2010,pu2021stochastic}. More recently, gradient tracking has emerged as a promising approach for distributed optimization \cite{Qu2017, nedic2017achieving, xi2018add}. For directed graphs, several gradient tracking methods have been proposed, including subgradient-push \cite{Tsianos2012, nedic2015distributed}, ADD-OPT \cite{xi2018add}, and Push-DIGing \cite{nedic2017achieving}, which rely on constructing a column-stochastic weight matrix that requires knowledge of agents' out-degree. Other algorithms, such as \cite{xi2018linear} and FROST \cite{Xin2019FROSTFastRO}, only use row-stochastic weights. However, these algorithms typically require separate iterations for Perron eigenvector estimation, which can cause stability issues. 

Recent works \cite{xin2018linear, pshi21, Saadatniaki2020} have addressed this issue with the algorithm $AB$-Push/Pull, which uses both row- and column-stochastic weights and removes the need for eigenvector estimation. Follow-up works, such as \cite{Saadatniaki2020,Angelia2022AB,Nguyen2023AccAB}, further establish linear convergence of this method for time-varying directed graphs. While references \cite{xin2018linear,Saadatniaki2020} derive linear convergence using arbitrary norms and the convergence bounds are not sharp, works by \cite{Angelia2022AB,Nguyen2023AccAB} provide an improved analysis using time-varying weighted average and weighted norms. The stochastic versions of the above-mentioned gradient tracking algorithms are studied in \cite{QureshiSADDOPT,XinSahuKhanKar2019}. In particular, S-ADDOPT is the stochastic version of ADD-OPT, while $\mathcal{S\!\!-\!\!AB}$ is the stochastic version of the $AB$-Push/Pull algorithm, considered over a static directed graph. Of significant relevance is \cite{XinSahuKhanKar2019}. Note that  \SABTV~ generalizes \cite{XinSahuKhanKar2019} by considering  general directed time-varying graphs and addressing the challenges posed by the time-varying nature of mixing matrices.

\textbf{Contribution.} This paper presents a stochastic gradient descent algorithm, \SABTV, for distributed optimization over directed time-varying communication graphs. Our algorithm is tailored to situations where the data is distributed across agents, each with access to a stochastic first-order oracle, and the gradient noise is assumed to have zero-mean and bounded variance. The \SABTV~algorithm leverages gradient tracking to mitigate the global variance resulting from the heterogeneity of the nodes' data. To facilitate directed communication networks, a two-phase update mechanism is employed, which involves pulling information about the decision variable from neighbors and pushing information about the gradients to neighbors using row- and column-stochastic weights, respectively. In addressing the challenge posed by the time-varying nature of the mixing matrices, our analysis incorporates time-varying weighted averages and norms associated with the row- and column-stochastic weight matrices. Using these norms, we establish consensus contractions for each update step and show that \SABTV~converges linearly to a neighborhood of the global minimizer in expectation for continuously-differentiable, strongly-convex, local cost functions, provided the constant step-size is sufficiently small.

\textbf{Notations.} 
All vectors are viewed as column vectors unless otherwise stated. The transpose of a vector $u\in \re^n$ is denoted by $u^\T$. The vector with all entries $1$ is represented $\one$. The $i$-th entry of a vector $u$ is denoted by $u_i$, while it is denoted by $[u_k]_i$ for a time-varying vector $u_k$. We denote $\min(u)=\min_i u_i$ and $\max(u)=\max_i u_i$. 
A nonnegative vector is called stochastic if its entries sum to $1$.

We use $[A_k]_{ij}$ to denote the $ij$-th entry of a matrix $A_k$. The notation $A \!\le\! B$ is used when $A_{ij} \!\le\! B_{ij}$ for all $i, j$, where $A$ and $B$ are matrices of the same dimension. A matrix $A$ is nonnegative if all its entries are nonnegative and ${\min}^{+}(A)$ denotes the smallest positive entry of $A$. A nonnegative matrix $A\in\mathbb{R}^{n\times n}$ is row-stochastic if $A\mathbf{1}=\mathbf{1}$, and a nonnegative matrix $B\in\mathbb{R}^{n\times n}$ is column-stochastic if $\mathbf{1}^{\T} B=\mathbf{1}^{\T}$. The identity matrix is denoted by $\mbi$. 

Given a positive vector $\ba=(a_1,\ldots,a_n)\in\re^n$, we denote:
\begin{center}
    $\la \bu,\bv\ra_{\ba}=\sum_{i=1}^m a_i\la u_i,v_i \ra $ and $\|\bu\|_{\ba}=\sqrt{\sum_{i=1}^m a_i\|u_i\|^2},$
\end{center}
where $\bu\!:=\![u_1,\ldots,u_n]^\T, \bv \!:=\![v_1,\ldots,v_n]^\T \!\!\in\!\! \re^{n\times p}$, and $u_i,v_i\!\in\!\re^p$. When $\ba = \one$, we write $\la \bu,\bv\ra$ and $\|\bu\|$. 

We let $[n]=\{1,\ldots,n\}$ for an integer $n\ge 1$. Given a directed graph $\G=([n],\mathcal{E})$, specified by the set of edges $\mathcal{E}\subseteq [n]\times[n]$ of ordered pairs of nodes, the in-neighbor and out-neighbor set for every agent $i$ are defined, as follows:
\[\cNini=\{j\in[n]|(j,i)\!\in\!\mathcal{E}\} \!\!\!\text{and}\!\!\! \cNouti=\{\ell\in[n]|(i,\ell)\!\in\!\mathcal{E}\}.\] 
A directed graph $\G$ is {\it strongly connected} if there is a directed path from any node to all other nodes in $\G$. 
We use $\mathsf{D}(\G)$ and $\mathsf{K}(\G)$ to denote the diameter and the maximal edge-utility of a strongly connected directed graph $\G$, respectively, as defined in \cite{Angelia2022AB,nguyen2022distributed}. We use a subscript $k$ to indicate the time instance.

\section{Stochastic $AB$/Push-Pull Method}
\label{sec:formu}
\subsection{Problem Formulation} 
\label{sec:Formulation}
Consider a system of $n$ agents connected over a communication network, with the aim of collaboratively solving the following optimization problem:
\begin{equation} \label{eq-problem}
\min_{x\in \re^p}~ f(x)=\frac{1}{n}\sum\limits_{i=1}^n f_i(x),
\end{equation}
where each function $f_i: \mathbb{R}^p \rightarrow \mathbb{R}$ represents the local cost function that is only known to agent $i$. We make the following assumption regarding these individual objective functions: 
\begin{assumption} \label{asm-functions}
Each $f_i$ is continuously differentiable and has $L$-Lipschitz continuous gradients, i.e., for some $L>0$,
\begin{equation*}
\|\nabla f_i(x)-\nabla f_i(y)\|\le L \|x-y\|,
\quad\hbox{for all $x,y\in \re^p$}.
\end{equation*}
\end{assumption}
	
\begin{assumption}\label{asm-strconv}
The average-sum function $f=\frac{1}{n}\sum_{i=1}^n f_i$ is $\mu$-strongly convex, i.e., for some $\mu>0$, 
\[\la \nabla f(x)-\nabla f(y),x-y\ra \ge \mu\|x-y\|^2\quad\hbox{for all $x,y\in \re^p$}.\]
\end{assumption}
 Assumption~\ref{asm-strconv} implies that problem~\eqref{eq-problem} has a unique optimal solution $x^*$ such that $x^*=\argmin_{x\in\re^p} f(x)$.

The agents aim to jointly find a globally optimal solution $x^*$ by performing local computations and exchanging information through a sequence of time-varying directed communication networks denoted by $\{{\bbG_k}\}$. At each time step $k$, agents communicate over a directed graph $\bbG_k=([n],\mathcal{E}_k)$, where each edge $(j,i)$ in $\mathcal{E}_k$ indicates that agent $i$ receives information from agent $j$. We consider the following assumption regarding the interaction graph of agents:
\begin{assumption} \label{asm-graphs}
For each $k$, the directed graph $\G_k$ is strongly connected and has a self-loop at every node $i\in[n]$.
\end{assumption}

Assumption~\ref{asm-graphs} can  be relaxed by considering a sequence of graphs that are $C$-strongly connected, i.e., for every $k\ge 0$, there exists an integer $C\ge 1$ such that the graph formed by the edge set $\mathcal{E}^C_k=\bigcup_{i=kC}^{(k+1)C-1}\mathcal{E}_i$ is strongly connected.

To solve Problem \eqref{eq-problem}, we assume that each agent $i$ is able to query a stochastic first-order oracle to obtain noisy gradient samples of the form $g_i(x,\xi_i)$ that satisfies the following standard assumptions specified in \cite{pu2018distributed}:
\begin{assumption} \label{asm-SFO}
For all $i\in[n]$ and $x\in\re^p$, each random vector $\xi_i\in\re^m$ is independent, and
\begin{flalign*}
&\text{(a)}\!\!\E[g_i(x,\xi_i)|x] = \nabla f_i(x),&&\\
&\text{(b)}\!\! \E[\|g_i(x,\xi_i)-\nabla f_i(x)\|^2|x] \le \sigma^2, \!\!\text{for some $\sigma>0$.}&&\!\!\!\!
\end{flalign*}
\end{assumption}
Assumption~\ref{asm-SFO} is widely applicable and particularly pertinent in several fields such as online distributed learning, reinforcement learning, generative models, and parameter estimation in signal processing and communication systems \cite{pu2018distributed}. 


\subsection{The \SABTV~Algorithm} \label{sec:algo}
This section introduces \SABTV, a distributed stochastic gradient tracking method proposed for large-scale systems where a centralized optimization approach is impractical. Specifically, each agent $i\in\{1,2,\ldots,n\}$ maintains a local copy $x_k^i\in \re^p$ and a direction $y_k^i\in\mathbb{R}^p$ which is an estimate of the ``global update direction" at iteration $k$. These variables are maintained and updated over time. Initially, each agent $i$ initializes their updates using arbitrary vectors $x_0^i$ and $y_0^i=g_i(x_0^i,\xi_0^i)$, without the need for coordination among agents. At time $k\ge 0$, every agent $i$ sends its vector $x_k^i$ and a scaled direction $[B_k]_{ji}y_k^i$ to its out-neighbors $j\in\Nikout$ and receives these vectors sent by its in-neighbors $j\in\Nikin$.
Upon the information exchange, every agent $i$ updates as follows:
\begin{subequations}\label{eq-met}
	\begin{align}
	&x_{k+1}^i =  \sum_{j=1}^n[A_k]_{ij}x_{k}^j  - \a y_k^i,\label{eq-x}\\
	&y_{k+1}^i =  \sum_{j=1}^n [B_k]_{ij}y_{k}^j + g_i(x_{k+1}^i,\xi_{k+1}^i) -g_i(x_k^i,\xi_k^i),\label{eq-y}
	\end{align}
\end{subequations}
where $\a>0$ is a constant step-size. The procedure is summarized in Algorithm 1. 

The updates for each agent are governed by two non-negative weight matrices $A_k$ and $B_k$ that  {\it align} with the topology of the graph $\bbG_k$, in the following sense:
\begin{align}\label{eq-alignA}
\!\!\![A_k]_{ij}\! >0,~\forall j\!\in\!\Nikin\!\cup\!\{i\}&;~\,[A_k]_{ij}\!=0,~\forall j\!\not\in\!\Nikin\!\cup\!\{i\},\!\!\\
\label{eq-alignB}
\!\!\![B_k]_{ji}\! >0,~\forall j\!\in\!\Nikout\!\cup\!\{i\}&;~\, [B_k]_{ji}\!=0,~\forall j\!\not\in\!\Nikout\!\cup\!\{i\}.\!\!
\end{align}
Each agent $i$ independently decides the entries $[A_k]_{ij}$ for their in-neighbors $j\in\Nikin$ while agent $j\in\Nikin$ determines the value $[B_k]_{ij}$. We make the following assumptions:

\begin{assumption} \label{asm-amatrices}
For each $k$, the matrix $A_k$ is row-stochastic, i.e., $A_k\1=\1$, and compatible with the graph $\bbG_k$ in the sense of relation~\eqref{eq-alignA}. Moreover, there exists a scalar $a\!>\!0$ such that $\min^+(A_k)\ge a$ for all $k\ge0$.
\end{assumption}

\begin{assumption}\label{asm-bmatrices}
For each $k$, the matrix $B_k$ is column-stochastic, i.e., ${\1}^{\T} B_k={\1}^{\T}$ and compatible with the graph $\bbG_k$ in the sense of relation~\eqref{eq-alignB}. Moreover, there exists a scalar $b>0$ such that $\min^+(B_k)\ge b$ for all $k\ge0$.
\end{assumption}


We can also write \eqref{eq-met} in the following compact form:
\begin{subequations}\label{eq-met-comp}
	\begin{align}
	&\bx_{k+1} =  A_k\bx_k  - \a \by_{k},\label{eq-x-comp}\\
	&\by_{k+1} =  B_k\by_{k} + \bg(\bx_{k+1},\bxi_{k+1}) -\bg(\bx_k,\bxi_{k}),\label{eq-y-comp}
	\end{align}
\end{subequations}
where $\bx_k\!=\![x_k^1,\ldots,x_k^n]^\T\!\in\!\re^{n\times p}$, $\by_{k}\!=\![y_k^1,\ldots,y_k^n]^\T\!\in\!\re^{n\times p}$ and $\bg(\bx_k,\bxi_{k})=[g_1(x_{k}^1,\xi_{k}^1),\ldots,g_n(x_{k}^n,\xi_{k}^n)]^\T\!\in\!\re^{n\times p}$ with $\bxi_k\!=\![\xi_k^1,\ldots,\xi_k^n]^\T\!\in\!\re^{n\times m}$. 

Finally, we denote by $\F_k$ the $\sigma$-algebra generated by the set of random vectors $\{\bxi_0,\ldots,\bxi_{k-1}\}$, and define by $\E[\cdot|\F_k]$ the conditional expectation given $\F_k$. 



\section{Convergence Analysis of \SABTV }\label{sec:conv_results}
In this section, we formalize the convergence analysis. Due to space limitations, we provide all the proofs of our results in our online technical report \footnote{Technical report, [Online] Available: https://tinyurl.com/SABTVreport}. We consider three critical error terms in terms of their expected values, namely, (i) the optimality gap $\E[\|\hat{x}_k -x^*\|^2]$, (ii) the consensus error $\E[\|\bx_k-\hat{\bx}_k\|_{\phi_k}^2]$, and (iii) the gradient tracking error $\E[S^2(\by_{k},\pi_k)]$, with $S(\by_{k}\!,\pi_k)$ defined as follows:
\begin{align}
S(\by_{k},\pi_k) =\sqrt{\sum_{i=1}^n[\pi_k]_i\Bigg\|\frac{y_k^i}{[\pi_k]_i} -\sum_{j=1}^n y_{k}^j\Bigg\|^2}.\label{eq-x-S-quants}
\end{align}
Here, $\{\phi_k\}$ is the sequence of stochastic vectors satisfying $\phi_{k+1}^{\T}A_k=\phi_k^{\T}$ and $\{\pi_k\}$ is the sequence of stochastic vectors defined in as $\pi_{k+1}=B_k\pi_k$, initialized with $\pi_0=\tfrac{1}{n}\1$. 

We define the constants $\kappa_k\ge 1$, $\varphi_k\ge 1$, $\gamma_k\in(0,1]$, $\psi_k>0$, $\tau_k\in(0,1)$, $c_k\in(0,1)$, $\nu_k>0$ and $\zeta_k>0$ that will be used in the analysis, as follows
\begin{align}\label{eq-const-all}
&\kappa_k\!=\!\sqrt{\tfrac{1}{\min(\pi_{k})}},~\varphi_k\!=\!\sqrt{\!\tfrac{1}{\min(\phi_{k})}},~\g_k\!=\!\!\sqrt{\max_{i\in[n]} ([\phi_{k+1}]_i[\pi_k]_i)},
 \nonumber\\
&\psi_k\!=\!n(\kappa_{k+1}^2\!\!-\!1),~ \tau_k \!=\!\sqrt{\!1 \!-\!   \tfrac{\min^2(\pi_k)\,b^2}
{\max^2(\pi_k\!) \max(\pi_{k+1}\!) \mathsf{D}(\bbG_k\!)\mathsf{K}(\bbG_k\!)}},\nonumber\\
& c_k=\sqrt{1- \tfrac{\min(\phi_{k+1})\, a^2} {\max^2(\phi_k)\,\mathsf{D}(\bbG_k)\mathsf{K}(\bbG_k)}},~\nu_k=\tfrac{4L^2\kappa_{k+1}^2\tau^2}{1-\tau^2}+2\psi_kL^2\!, \nonumber\\
&\zeta_k = \left(c\varphi_{k+1}+\varphi_k \right)^2\nu_k,
\end{align}
where $c\in(0,1)$ and $\tau\in(0,1)$ are upper bounds for $c_k$ and $\tau_k$, respectively. Furthermore, 
we assume that $\eta$ is the lower bound for $\phi_{k+1}^\T\pi_k$, thus, $1<\eta<2$. We let $\psi>0$, $\kappa>1$, and $\varphi>1$, be upper bounds for $\psi_k$, $\kappa_k$, and $\varphi_k$, respectively. 

Let Assumption~\ref{asm-functions}, Assumption~\ref{asm-strconv}, Assumption~\ref{asm-graphs}, Assumption~\ref{asm-SFO}, Assumption~\ref{asm-amatrices} and Assumption~\ref{asm-bmatrices} hold. In the following propositions, we establish bounds on the conditional expectations of $\|\hat{x}_k -x^*\|^2$, $\|\bx_k-\hat{\bx}_k\|_{\phi_k}^2$ and $S^2(\by_{k},\pi_k)$, respectively.

\begin{table}[t!]
\vspace{0.2cm}
\centering \normalsize
    \begin{tabular}{l}
    \hline
    \multicolumn{1}{c}{\textbf{Algorithm 1: \SABTV}}\\
    \hline
    Every agent $i\in[n]$ initializes with arbitrary initial vectors\\ $x_0^i\in\re^{p}$ and $y_0^i=g_i(x_0^i,\xi_0^i)$. \\
    Agents are instructed to use the step-size $\a$.\\
    \textbf{for} $k=0,1,\ldots,$ every agent $i\in[n]$ does the following:\!\!\\
    \emph{  } Receives $x_{k}^j$ and $[B_k]_{ij}y_k^j$ from in-neighbors $j\in\cNinik$;\\
    \emph{  } Sends $x_k^i$ and $[B_k]_{ji}y_k^i$ to  out-neighbors $j\in\cNoutik$;\\
    \emph{ } Chooses the weights $[A_k]_{ij},j\in\cNinik$;\\
    \emph{ } Updates the action $x_{k+1}^i$ using \eqref{eq-x};\\
    \emph{ } Updates the direction $y_{k+1}^i$ using \eqref{eq-y};\\
    \textbf{end for}\\
    \hline
    \end{tabular}
    \vspace{-0.3cm}
\end{table}

\begin{proposition} \label{prop-waverx}
For $0\!<\!\a\!<\!\tfrac{2}{n\eta(L+\mu)}$,
we have for all $k\ge0$, 
\begin{align*}
&\E[\|\hat x_{k+1}-x^*\|^2|\F_k]\nonumber\\
\le& \!\left(\!1\!-\! \tfrac{\a n\mu \phi_{k+1}^\T\pi_k}{2}\right)\!\|\hat x_k - x^*\|^2\!+\!\tfrac{3\a \phi_{k+1}^\T\pi_kL^2\varphi_k^2}{\mu}\|\bx_k-\hat{\bx}_k\|_{\phi_k}^2\nonumber\\
&+\tfrac{3\a}{ n \mu \phi_{k+1}^\T\pi_k}\E[S^2(\by_{k},\pi_k)|\F_k]
+\tfrac{3\a^2 n (\phi_{k+1}^\T\pi_k)^2\sigma^2}{2}.
\end{align*}
\end{proposition}

\begin{proposition}\label{prop-xcontract}
We have for all $k\ge0$,
\begin{align*}
\E\left[\left\|\bx_{k+1}-\hat \bx_{k+1}\right\|_{\phi_{k+1}}^2|\F_k\right]
\le \tfrac{2nL^2\varphi_k^2(1+c^2)\g_k^2\a^2}{1-c^2} \|\hat x_k - x^*\|^2\\
\left(\tfrac{1+c^2}{2}+\tfrac{2nL^2\varphi_k^2(1+c^2)\g_k^2\a^2}{1-c^2}\right)\left\|\bx_k - \hat \bx_k\right\|_{\phi_{k}}^2\\
+\tfrac{\a^2 (1+c^2)\g_k^2}{1-c^2}\E[S^2(\by_{k},\pi_k)|\F_k]+\tfrac{2\a^2(1+c^2)\g_k^2\sigma^2}{1-c^2}.
\end{align*}
\end{proposition}

\begin{proposition}\label{prop-ycontract}
We have for all $k\ge0$,
\begin{align*}
\E\!\left[S^2(\by_{k+1}\!,\pi_{k+1})|\F_k\right]
\le (\zeta_k+2nL^2\varphi_k^2\nu_k \a^2)\|\bx_k-\hat{\bx}_k\|_{\phi_k}^2\! \\
\!+2nL^2\varphi_k^2\nu_k \a^2 \|\hat x_k - x^*\|^2\!+\!(4n\psi_k\!+\!2\a Ln\psi_k\!+\!2\nu_k \a^2)\sigma^2\nonumber\\
+\left(\tfrac{1+\tau^2}{2}+\nu_k \a^2\right)\E[S^2(\by_{k},\pi_k)|\F_k].
\end{align*}
\end{proposition}

\subsection{Composite Relation} \label{subsec-comrel}
Defining the vector $V_k$ as follows:
\begin{align*}\label{eq-vk}
\!\!V_k=\Big(\E[\|\hat{x}_k -x^*\|^2],\E[\|\bx_k-\hat{\bx}_k\|_{\phi_k}^2],\E[S^2(\by_{k},\pi_k)]\Big)^{\!\T}\!\!\!.
\end{align*}
Then, $V_k$ follows the dynamical system below:
\begin{equation}\label{eq-vkrel}
V_{k+1}\le M(\a)V_k+\bb(\a), \qquad \hbox{for all $k\ge0$},
\end{equation}
where
\begin{align}\label{eq-gmatrixm}
M(\a)&=\!\left[\begin{array}{ccc}
\!1-m_{1}\a & m_{2}\a & m_{3}\a \cr
m_{4}\a^2 &\!\tfrac{1+c^2}{2}+ m_{4}\a^2\! & m_{5}\a^2\cr
m_{6}\a^2 & m_{7}+m_{6}\a^2 \!& \!\tfrac{1+\tau^2}{2}+m_{8}\a^2 \end{array}\right]\!,\!\!\\
\!\!\!\!\!\!\text{and} &\bb(\a)=\left[b_{1}\a^2,b_{2}\a^2,b_{3}+b_{4}\a+b_{5}\a^2\right]^\T. \label{eq-vecg}
\end{align} 

The constants in $M(\a)$ and $\bb(\a)$ are given by
\begin{align*}
&m_{1}=\tfrac{n\mu}{2},~m_{2}=\tfrac{3L^2\varphi^2}{\mu},~m_{3}=\tfrac{3}{ n \mu \eta},~m_{4}=\tfrac{2nL^2\varphi^2(1+c^2)}{1-c^2},\\
&m_{5}=\tfrac{1+c^2}{1-c^2},~m_{6}=2nL^2\varphi^2\nu,~m_{7}=\zeta,~m_{8}=\nu,~b_{1}=\tfrac{3n\sigma^2}{2}\!\!,\\
&b_{2}=\tfrac{2(1+c^2)}{1-c^2},~b_{3}=4n\psi\sigma^2,~b_{4}=2Ln\psi\sigma^2,~b_{5}=2\nu\sigma^2,
\end{align*}
with $\nu=\tfrac{4L^2\kappa^2\tau^2}{1-\tau^2}+2\psi L^2$ and $\zeta = 4\varphi^2\nu$.

\subsection{Convergence Theorem} \label{subsec-conv-theo}
We now present the main result of this paper, which states that the \SABTV~algorithm (Algorithm 1) converges linearly to a neighborhood of the global minimizer. 
\begin{theorem}\label{theo-main}
Let Assumption~\ref{asm-functions}, Assumption~\ref{asm-strconv}, Assumption~\ref{asm-graphs}, Assumption~\ref{asm-SFO}, Assumption~\ref{asm-amatrices} and Assumption~\ref{asm-bmatrices} hold. Consider the iterates produced by Algorithm 1, the notations in \eqref{eq-const-all}, the matrix $M(\a)$ as defined in \eqref{eq-gmatrixm} and the vector $\bb(\a)$ as in \eqref{eq-vecg}. If the step-size $\a>0$ is chosen such that 
\begin{align}\label{eq:alpha-range}
\!\!\a \le \min \!\Big\{\!&\tfrac{2}{n\eta(L+\mu)},\!\tfrac{\mu\eta(1-c^2)(1-\tau^2)}{4L\varphi\sqrt{\!(n\mu^2\eta^2+24L^2\varphi^2)(\eta+16\kappa^2)\!}},\nonumber\\
 &\tfrac{\mu\eta(1-\tau^2)}{\kappa L\sqrt{2n\mu^2\eta^2(n+4)+6L^2\varphi^2(n\eta+32)}}\Big\},
\end{align}
then $\rho_{M}<1$ where $\rho_{M}$ is the spectral radius of $M(\a)$  and, the vector $(\mbi-M(\a))^{-1}\bb(\a)$ has non-negative components, and we have that 
\begin{align*}
&\lim\sup_{k\to\infty} V_k \le (\mbi-M(\a))^{-1}\bb(\a),
\end{align*}
with a linear convergence rate of the order of $\mathcal{O}\Big(\rho_M^k\Big)$.
\end{theorem}
\allowdisplaybreaks
\begin{proof}
Recall that from \eqref{eq-vkrel} we have
\[V_{k+1}\le M(\a)V_k+\bb(\a), \qquad \hbox{for all $k\ge0$}.\]
The goal is to find the range of $\a$ such that the spectral radius $\rho_{M}$ of $M(\a)$ satisfies $\rho_{M}<1$. It suffices to solve for the range of $\a$ such that $M(\a)\bdelta < \bdelta$ holds for some positive vector $\bdelta=[\delta_1,\delta_2,\delta_3]^\T$. Expanding $M(\a)\bdelta < \bdelta$, we obtain
\begin{align}
(-m_1\delta_1+m_2\delta_2+m_3\delta_3)\a~&<0,\label{eq-alpha1}\\
(m_{4}\delta_1+m_{4}\delta_2+m_{5}\delta_3)\a^2&<\tfrac{1-c^2}{2}\delta_2,\label{eq-alpha2}\\
(m_{6}\delta_1+m_{6}\delta_2+m_{8}\delta_3)\a^2&<\tfrac{1-\tau^2}{2}\delta_3-m_{7}\delta_2.\label{eq-alpha3}
\end{align}
It is straightforward to check that the following $\delta_1$, $\delta_2$ and $\delta_3$ satisfy \eqref{eq-alpha1} and and make the right-hand side of \eqref{eq-alpha3} positive:
\begin{align*}
\delta_1=\tfrac{2}{m_1}\left(m_2+\tfrac{4m_3m_7}{1-\tau^2}\right),~ \delta_2=1,~ \delta_3=\tfrac{4m_7}{1-\tau^2}.
\end{align*}
Plugging in these values of $\delta_1,\delta_2,\delta_3$ into \eqref{eq-alpha2} and \eqref{eq-alpha3}, 
we have $M(\a)\bdelta < \bdelta$ when the following inequalities hold
\begin{align}
&\!\!\left(\!m_4+\tfrac{2m_4}{m_1}\!\left(\!m_2+\tfrac{4m_3m_7}{1-\tau^2}\!\right)\!+\tfrac{4m_5m_7}{1-\tau^2}\!\right)\!\a^2<\tfrac{1-c^2}{2}\!,\!\!\!\label{eq-solve-alpha-range1}\\
&\!\!\left(\!\tfrac{m_6}{m_7}+\tfrac{2m_6}{m_1m_7}\!\left(\!m_2+\tfrac{4m_3m_7}{1-\tau^2}\!\right)\!+\tfrac{4m_8}{1-\tau^2}\!\right)\!\a^2<1.\label{eq-solve-alpha-range2}
\end{align}
We can further verify that 
\begin{align*}
m_4+\tfrac{2m_4}{m_1}&\left(m_2+\tfrac{4m_3m_7}{1-\tau^2}\right)\!+\tfrac{4m_5m_7}{1-\tau^2}\nonumber\\
&<\tfrac{8L^2\varphi^2(n\mu^2\eta^2+24L^2\varphi^2)(\eta+16\kappa^2)}{\mu^2\eta^2(1-c^2)(1-\tau^2)^2},\\
\tfrac{m_6}{m_7}+\tfrac{2m_6}{m_1m_7}&\!\left(m_2+\tfrac{4m_3m_7}{1-\tau^2}\right)+\tfrac{4m_8}{1-\tau^2}\nonumber\\
&< \tfrac{\kappa^2L^2\left[2n\mu^2\eta^2(n+4)+6L^2\varphi^2(n\eta+32)\right]}{\mu^2\eta^2(1-\tau^2)^2}.
\end{align*}
Using the above bounds in \eqref{eq-solve-alpha-range1} and \eqref{eq-solve-alpha-range2}, we obtain the upper bounds for $\a$ as in \eqref{eq:alpha-range} which establishes the proof.
\end{proof}
\begin{remark}
The error bounds in Theorem~\ref{theo-main} go to zero as the step-size $\a$ gets smaller and the variance $\sigma$ on the gradient noise decreases.
\end{remark}

\section{Numerical results} \label{sec:simulation}
In this section, we examine a binary classification problem. We consider $n=10$ agents embedded in a communication network that is directed and time-varying, and each agent represents a node on the graph. 
To ensure the graphs are strongly connected, we establish a directed cycle linking all agents at each iteration.

We consider a total of $N$ labeled data points for training, with each node $i$ possessing a local batch of $m_i$ training samples. The $j$-th sample at node $i$ is a tuple $\{b_{ij},y_{ij}\} \!\subseteq \!\re^p \!\times \!\{
+1,\!-1\}$. To construct an estimate of the coefficients $x\!\!=\!\![x_0,x_{1:}^\T]^\T\!\!\in\! \re^{p+1}$, where $x_{1:}\!\!=\!\![x_1,\ldots,x_p]^\T\!$, we will use the principle of maximum likelihood and define the local logistic regression cost function $f_{ij}$ for the $j$-\textit{th} image at node $i$ as:
\[f_{ij}(x) = \ln \Big[1+\exp\Big\{-(x_{1:}^\T b_{ij}+x_0)y_{ij}\Big\} \Big] \!+\dfrac{\lambda}{2}\|x\|^2\!,\]
which is smooth and strongly convex due to the inclusion of the L2-regularization. The agents cooperatively solve the following optimization problem:
\[\min_{x\in \mathbb{R}^{p+1}}~ \frac{1}{n}\sum\limits_{i=1}^n \dfrac{1}{m_i}\sum_{j=1}^{m_i}f_{ij}(x).\]

We perform digit classification on the widely used MNIST dataset, with the objective of accurately distinguishing between handwritten digits $3$ and $7$. 
We use $12000$ samples for training and $2400$ samples for testing purposes. 


We evaluate the performance of centralized and distributed methods. The centralized gradient descent (CGD) method employs the entire batch, computing $12000$ gradients at each iteration. In contrast, the centralized stochastic gradient descent (CSGD) approach uses a single data point, randomly sampled from the entire batch, for each iteration. We also evaluate the performance of the non-stochastic $AB$/Push-Pull (ABPP) method for distributed algorithms \cite{Saadatniaki2020, Angelia2022AB}, where each agent processes a local batch of $1200$ labeled data points. On the other hand, 
in the proposed  \SABTV~algorithm,
each agent selects one data point uniformly from its local batch. 

\begin{figure}[ht!]
\vspace{-0.3cm}
	\centering
	\hspace*{-1.8em}
		\subfigure{
	 \includegraphics[width=0.23\textwidth]{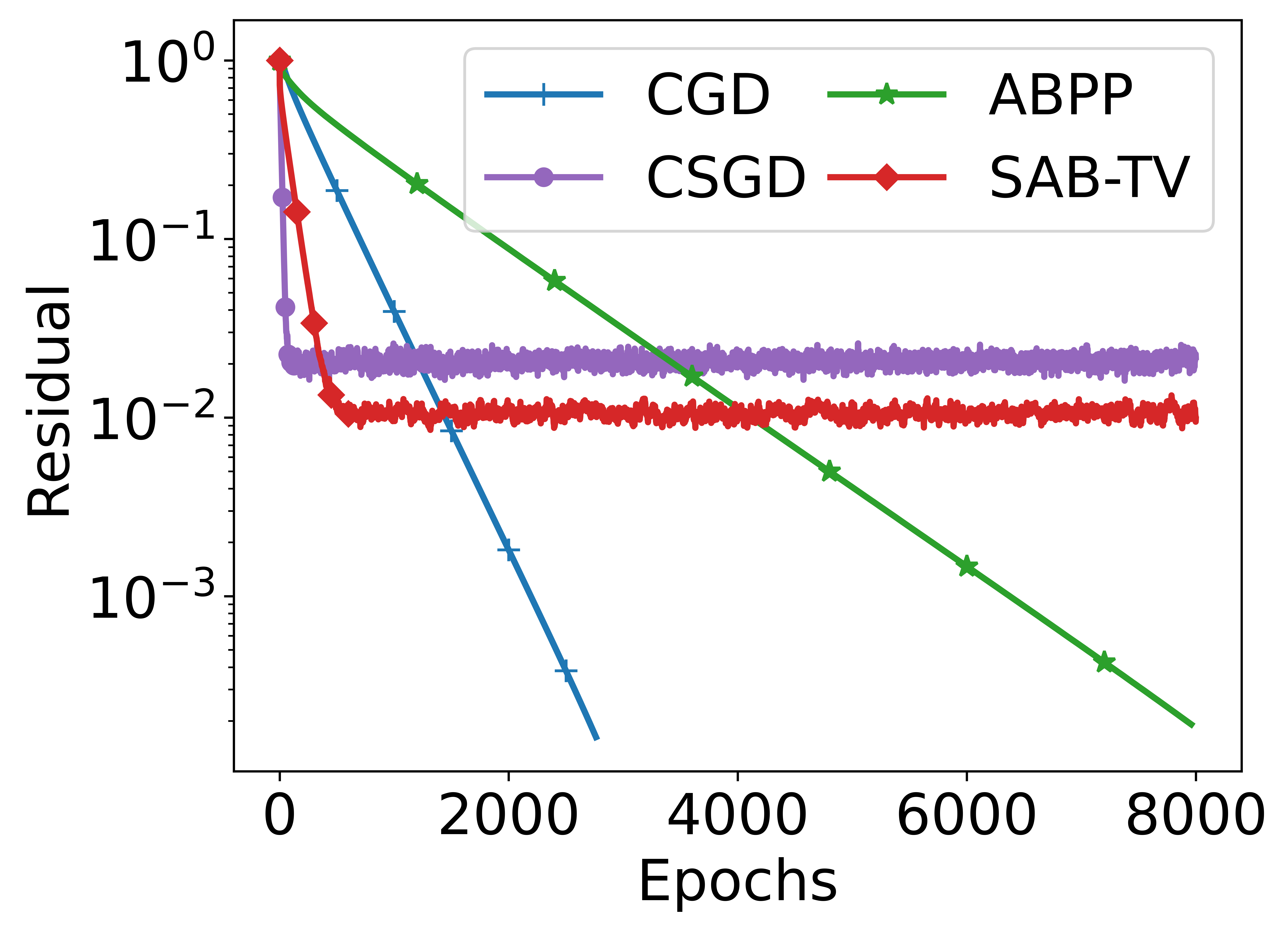}
	    \label{fig:Residual}  } 
	   \hspace*{-1em} 
		 \subfigure{
	\includegraphics[width=0.23\textwidth]{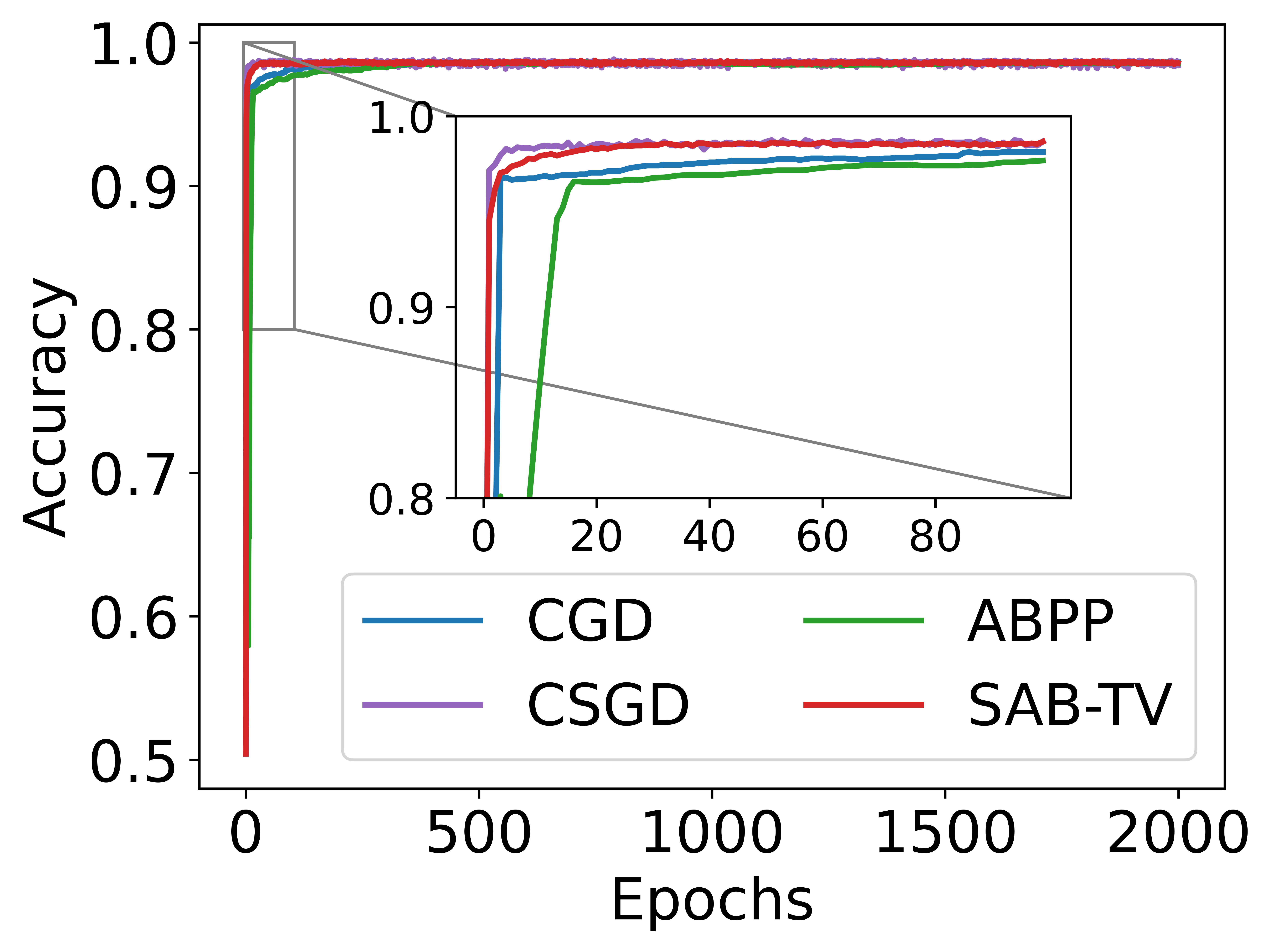}
	   \label{fig:Accuracy}
	} \hspace*{-1.8em}
 \vspace{-0.3cm}
	\caption{Residual (left) and test accuracy (right).}\label{fig:ResidualAccuracy}
 \vspace{-0.2cm}
\end{figure}

The residuals are displayed in Figure~\ref{fig:Residual}, while  Figure~\ref{fig:Accuracy} compares the accuracy of all four algorithms on the test set.
The $x$-axis represents the number of epochs,
with each epoch corresponding to computations on the entire batch. We observe that \SABTV~outperforms $AB$-Push/Pull initially, which aligns with the performance of their centralized counterparts, CSGD and CGD. However, due to the inexact gradient, both \SABTV~and CSGD  only converge to a vicinity of the optimal solution, whereas the $AB$-Push/Pull and CGD algorithms exhibit linear convergence to the optimal solution. After 50 epochs, all algorithms achieve a test set accuracy of over $97\%$.

\section{Conclusions}
\label{sec:conc}
This paper introduces a stochastic gradient-based algorithm for distributed optimization over directed time-varying graphs where data is distributed across multiple agents, each having access to only inexact gradients of its local cost function. The \SABTV~algorithm converges linearly to a neighborhood of the optimum in expectation, for a sufficiently small constant step-size when the local cost functions are smooth and strongly convex. We also conducted numerical simulations based on the MNIST dataset to 
visualize the theoretical results. 

\bibliographystyle{IEEEtran}
\bibliography{references}

\end{document}